\documentclass[12]{article}
\usepackage{fleqn}
\usepackage{graphicx}
\usepackage{amssymb}
\begin{document}
\Large
\begin{center}
{\bf A Classification of the Veldkamp Lines of the\\ Near Hexagon L$_3 \times$
GQ(2,\,2)}
\end{center}
\large
\vspace*{-.2cm}
\begin{center}
R.M. Green$^{1}$ and Metod Saniga$^{2}$
\end{center}
\vspace*{-.5cm} \normalsize
\begin{center}
$^{1}$Department of Mathematics, University of Colorado Boulder, Campus
Box 395\\ Boulder CO 80309-0395,
U. S. A. \\
(rmg@euclid.colorado.edu)

\vspace{0.3cm}

$^{2}$Astronomical Institute, Slovak Academy of Sciences\\
SK-05960 Tatransk\' a Lomnica, Slovak Republic\\
(msaniga@astro.sk)

\vspace*{-.2cm}


\end{center}

\vspace*{.0cm} \noindent \hrulefill

\vspace*{.0cm} \noindent {\bf Abstract}\\
Using a standard technique sometimes (inaccurately) known as Burnside's Lemma, it is shown that the Veldkamp space of the near hexagon $L_3 \times$GQ(2,\,2) features 156 different types of lines.
We also give an explicit description of each type of a line by listing the types of the three geometric hyperplanes it consists of and describing the properties of its core set, that is the subset of points of $L_3 \times$GQ(2,\,2) shared by the three geometric hyperplanes in question.

\vspace*{-.2cm} \noindent \hrulefill

\vspace*{.1cm}
\section{Introduction}

Brouwer {\it et al.} \cite{bchw} proved that there are eleven isomorphism types
of slim dense near hexagons.  Of these eleven, the near hexagons of sizes $27$,
$45$ and $81$ are the most promising for physical applications.  This paper
is devoted to a study of the second of these three examples and its Veldkamp
space.  The first of the three examples was described in our paper \cite{gs}, and 
we plan to study the third case in a future work.
The $45$ point space we study here is the product $L_3 \times$GQ(2,\,2), 
where $L_3$ is the line containing three points and GQ(2,\,2) is the generalized quadrangle of order two.

\section{Near polygons, quads, geometric hyperplanes and\\ Veldkamp spaces}
In this section we gather all the basic  notions and well-established theoretical results that will be
needed in the sequel.

A {\it near polygon} (see, e.\,g., \cite{bruyn} and references
therein) is a connected partial linear space $S = (P, L, I)$, $I
\subset P \times L$, with the property that given a point $x$ and
a line $L$, there always exists a unique point on $L$ nearest to
$x$. (Here distances are measured in the point graph, or
collinearity graph of the geometry.)  If the maximal distance
between two points of $S$ is equal to $d$, then the near polygon
is called a near $2d$-gon. A near 0-gon is a point and a near
2-gon is a line; the class of near quadrangles coincides with the
class of generalized quadrangles.

A nonempty set $X$ of points in a near polygon $S = (P, L, I)$ is
called a subspace if every line meeting $X$ in at least two points
is completely contained in $X$. A subspace $X$ is called
geodetically closed if every point on a shortest path between two
points of $X$ is contained in $X$. Given a subspace $X$, one can
define a sub-geometry $S_X$ of $S$ by considering only those
points and lines of $S$ that are completely contained in $X$. If
$X$ is geodetically closed, then $S_X$ clearly is a
sub-near-polygon of $S$. If a geodetically closed sub-near-polygon
$S_X$ is a non-degenerate generalized quadrangle, then $X$ (and
often also $S_X$) is called a {\it quad}.

A near polygon is said to have order $(s, t)$ if every line is
incident with precisely $s+1$ points and if every point is on
precisely $t+1$ lines. If $s = t$, then the near polygon is said
to have order $s$. A near polygon is called {\it dense} if every
line is incident with at least three points and if every two
points at distance two have at least two common neighbours. A near
polygon is called {\it slim} if every line is incident with
precisely three points. It is well known (see, e.\,g.,
\cite{pay-thas}) that there are, up to isomorphism, three slim
non-degenerate generalized quadrangles. The $(3 \times 3)$-grid is
the unique generalized quadrangle of order $(2, 1)$, GQ$(2, 1)$.
The unique generalized quadrangle of order 2, GQ$(2, 2)$, is the
generalized quadrangle of the points and lines of PG(3, 2) that
are totally isotropic with respect to a given symplectic form.
The points and lines lying on a given nonsingular elliptic quadric
of PG$(5, 2)$ define the unique generalized quadrangle of order
$(2, 4)$, GQ$(2,4)$. Any {\it slim dense} near polygon contains
quads, which are necessarily isomorphic to either GQ$(2, 1)$,
GQ$(2, 2)$ or GQ$(2, 4)$.

Next, a {\it geometric hyperplane} of a partial linear space is a
proper subspace meeting each line (necessarily in a unique point
or the whole line). The set of points at non-maximal distance from
a given point $x$ of a dense near polygon $S$ is a hyperplane of
$S$, usually called the {\it singular} hyperplane (or {\it perp-set}) with {\it
deepest} point $x$. Given a hyperplane $H$ (or any subset of points $\mathcal{C}$) of $S$, one defines the
{\it order} of any of its points as the number of lines through
the point that are fully contained in $H$ ($\mathcal{C}$); a point of $H$ ($\mathcal{C}$) is called {\it deep} if all the lines
passing through it are fully contained in $H$ ($\mathcal{C}$). If $H$ is a hyperplane of a dense
near polygon $S$ and if $Q$ is a quad of $S$, then precisely one
of the following possibilities occurs: (1) $Q \subseteq H$; (2) $Q
\cap H = x^{\perp} \cap Q$ for some point $x$ of $Q$; (3) $Q \cap
H$ is a sub-quadrangle of $Q$; and (4) $Q \cap H$ is an ovoid of
$Q$. If case (1), case (2), case (3), or case (4) occurs, then $Q$
is called, respectively, {\it deep}, {\it singular}, {\it
sub-quadrangular}, or {\it ovoidal} with respect to $H$. If $S$ is
slim and $H_1$ and $H_2$ are its two distinct hyperplanes, then
the complement of symmetric difference of $H_1$ and $H_2$,
$\overline{H_1 \Delta H_2}$, is again a hyperplane; this means
that the totality of hyperplanes of a slim near polygon form a
vector space over the Galois field with two elements, ${\mathbb F}_2$. In
what follows, we shall put $\overline{H_1 \Delta H_2} \equiv H_1
\oplus H_2$ and call it the (Veldkamp) sum of the two hyperplanes.

Finally, we shall introduce the notion of the {\it Veldkamp space},
$\mathcal{V}(\Gamma)$, 
of a point-line incidence geometry $\Gamma(P, L)$ \cite{buek}.  
Here,  $\mathcal{V}(\Gamma)$  is the
space in  which (i) a point is a geometric hyperplane of  $\Gamma$
and (ii) a line is the collection $H'H''$ of all geometric
hyperplanes $H$ of $\Gamma$  such that $H' \cap H'' = H' \cap H =
H'' \cap H$ or $H = H', H''$, where $H'$ and $H''$ are distinct
points of $\mathcal{V}(\Gamma)$.
Following \cite{spph,sglpv}, we
adopt also here the definition of Veldkamp space given by
Buekenhout and Cohen \cite{buek} instead of that of Shult
\cite{shult}, as the latter is much too restrictive by requiring
any three distinct hyperplanes $H'$, $H''$ and $H'''$ of $\Gamma$
to satisfy the following two conditions: i) $H'$ is not properly contained
in $H''$ and ii) $H' \cap H'' \subseteq H'''$
implies $H' \subset H'''$ or  $H' \cap H'' = H' \cap H'''$. The two
definitions differ in the crucial fact that whereas the Veldkamp space in the sense
of Shult is {\it always} a linear space, that of Buekenhout and Cohen needs
not be so; in other words, Shult's Veldkamp lines are always
of the form $\{H \in \mathcal{V}(\Gamma) ~|~ H \supseteq H' \cap H''\}$ for certain geometric hyperplanes
$H'$ and $H''$.

\section{The near hexagon L$_3 \times$GQ(2,\,2)}
The near hexagon $L_3 \times$ GQ(2,\,2) has recently \cite{sng} caught an attention of theoretical physicists due to the fact that its main constituent, the generalized quadrangle GQ(2,\,2), 
reproduces the commutation relations of the 15 elements of the two-qubit 
Pauli group (see, e.\,g., \cite{ps}), with each of its ten embedded copies of GQ(2,\,1) playing, remarkably, the role of the so-called {\it Mermin magic square} \cite{merm} --- the smallest configuration of two-qubit observables furnishing a very important proof of contextuality of quantum mechanics.
A well-known construction of GQ(2,\,2) identifies the points with two-element
subsets of $\{1, 2, 3, 4, 5, 6\}$, with two points being collinear if and
only if they are equal or disjoint.
The natural action of $S_6$ on this set of size 6 induces automorphisms of
GQ(2,\,2).  In fact, when considered in this way, $S_6$ turns out to be the 
full automorphism group.

It is known that every geometric hyperplane of a slim dense near polygon arises
from its universal embedding.  It can be shown from this that,
equipped with the operation of Veldkamp sum, the Veldkamp space $V_{GQ(2,2)}$ 
is isomorphic to PG(4,\,2), the projective space obtained from a 5-dimensional 
space over ${\mathbb F}_2$ (see also \cite{spph}).  It follows that GQ(2,\,2) has $2^5 - 1 = 31$ 
geometric hyperplanes, which turn out to be of three types:

\begin{itemize}
\item[(i)] 15 perp-sets, with 7 points each;
\item[(ii)] 10 grids (copies of GQ(2,\,1)), with 9 points each;
\item[(iii)] 6 ovoids, with 5 points each.
\end{itemize}

In other words, there are three orbits of geometric hyperplanes under the 
action of $S_6$.

Identifying the points of GQ(2,\,2) with two-element subsets of the set
$\{1, 2, 3, 4, 5, 6\}$ as described earlier, we find that an example of
an ovoid is the set $$
e_1 := \{ 
\{1, 2\}, 
\{1, 3\}, 
\{1, 4\}, 
\{1, 5\}, 
\{1, 6\} \}
.$$  The other ovoids, $e_2, e_3, \ldots, e_6$ are obtained from $e_1$ by 
acting by the transposition $(1, i)$ for $i = 2, 3, \ldots, 6$ respectively.

The Veldkamp 
sum $e_i + e_j$ (for $1 \leq i < j \leq 6$) is the perp-set of the point
$\{i, j\}$.  If we have $$
\{1, 2, 3, 4, 5, 6\} = \{i, j, k, l, m, n\}
$$ in some order, then the sum $e_i + e_j + e_k$ is the grid whose elements are
the nine points $$
\{\{a, b\} : a \in \{i, j, k\} \mbox{\ and\ } b \in \{l, m, n\}\}
.$$  It follows that the six ovoids are a spanning set for $V_{GQ(2,2)}$.
Since each point of GQ(2,\,2) lies in precisely two ovoids, it follows that
we have the relation $$
e_1 + e_2 + e_3 + e_4 + e_5 + e_6 = 0
,$$ where $0$ denotes the subset of GQ(2,\,2) consisting of all $15$ points.
Since we have an isomorphism $V_{GQ(2,2)} \cong$ PG(4,\,2), it
follows by a counting argument that this is the only nontrivial dependence 
relation between the $e_i$, and thus that the ovoids $e_1, \ldots, e_5$ 
form a basis for $V_{GQ(2,2)}$.

The points of the near hexagon $L_3 \times$ GQ(2,\,2) are simply the 45 
ordered pairs $(p, q)$ where $p$ is a point of $L_3$ and $q$ is a point of 
GQ(2,\,2). 
We call a collection of 15 points $(p, q)$ sharing the same value of $p$
a {\it layer} of the near hexagon. A layer is an example of a quad in the sense
of \S2. We imagine that the points of $L_3$ are arranged vertically, and
we will sometimes use terms like ``the top quad'' to refer to 
one of the layers of the near hexagon.

Two points $(p_1, q_1)$ and $(p_2, q_2)$ of $L_3 \times$ GQ(2,\,2) are
collinear if either 
\begin{itemize}
\item[(i)] $p_1 = p_2$ and $q_1$ is collinear to $q_2$, or
\item[(ii)] $p_1$ is collinear to $p_2$ and $q_1 = q_2$.
\end{itemize}
The lines of $L_3 \times$ GQ(2,\,2) are of two types. The {\it type-one} lines are the
15 lines of the form $\{(p, q) : p \in L_3\}$ for a fixed point 
$q \in$ GQ(2,\,2).
The {\it type-two} lines
are the 45 lines of the form $\{(p, q) : q \in L\}$ for a fixed
$p \in L_3$ and some line $L$ of GQ(2,\,2). 

The near hexagon $L_3 \times$ GQ(2,\,2) has a number of obvious automorphisms.
One type of automorphism involves permuting the three GQ(2,\,2)-quads, but making no other changes.  The subgroup of all such
automorphisms is isomorphic to $S_3$.
Another type of automorphism involves acting diagonally on the three GQ(2,\,2)-quads
by $S_6$, the automorphism group of GQ(2,\,2).  This action commutes with
the action of $S_3$ just mentioned, and produces a group of automorphisms
isomorphic to $S_6 \times S_3$.
It turns out that this is the full automorphism group, as shown by
Brouwer {\it et al.} \cite{bchw}.

From now on, let us denote the Veldkamp space of $L_3 \times$ GQ(2,\,2)
by $V$.  Some features of $V$ are close to obvious, which stems from Sec.\,2. 
One of these is that the intersection of one of the three GQ(2,\,2)-quads with a 
point of $V$ (regarded as a subset of the 45 points) can take one of two 
forms.  Either the GQ(2,\,2)-quad is completely filled in (i.\,e., it is deep), or 
takes the form of one of the geometric hyperplanes of GQ(2,\,2) (i.\,e., it is singular, sub-quadrangular or ovoidal).
Furthermore, the Veldkamp sum of any two of the layers (regarded as subsets
of GQ(2,\,2) under some obvious identification) must be equal to the third layer.
It follows from this that $V$ contains $2^{10} - 1 = 1023$ points.  

The above discussion shows that, as an $S_6 \times S_3$-module over
${\mathbb F}_2$, $V$ is isomorphic to $M \otimes N$, where $M$ is the
$5$-dimensional module for $S_6$ described earlier, and $N$ is the 
$S_3$-module obtained by quotienting the $3$-dimensional permutation
module $\{f_1, f_2, f_3\}$ for $S_3$ by the submodule spanned by 
$f_1 + f_2 + f_3$.  The set $\{f_1, f_2\}$ then form a basis for $N$, and
the set $$
\{e_i \otimes f_j: 1 \leq i \leq 5, \ 1 \leq j \leq 2\}
$$ forms a basis for $V$.  We will write this basis for short as
$\{e_1, \ldots e_{10}\}$, where for $1 \leq i \leq 5$, $e_i$ denotes
$e_i \otimes f_1$, and for $6 \leq i \leq 10$, $e_i$ denotes
$e_{i-5} \otimes f_2$.

\section{The classification of hyperplanes}

The geometric hyperplanes of $L_3 \times$ GQ(2,\,2) were classified 
in \cite{sng}. Up to automorphisms, there are eight types of them,
denoted by $H_1$ to $H_8$ and described in detail in \cite[Table 2]{sng}.
We now explain how these eight types can be reconstructed using the results
in the previous section.

The description of the hyperplanes of GQ(2,\,2) above
can be used to identify each hyperplane with one of the 31 nontrivial
set partitions of a 6-element into two pieces. If $S$ and $T$ are disjoint
nonempty sets for which $$
S \cup T = \{1, 2, 3, 4, 5, 6\},
$$ then we identify the pair $\{S, T\}$ with the hyperplane $$
\sum_{i \in S} e_i = \sum_{j \in T} e_j
.$$ If $|S| \geq |T|$, we associate the
partition $(|S|, |T|)$ of the number $6$ to the set partition $\{S, T\}$.
Under these identifications, the partitions of $6$ given by $(5, 1)$,
$(4, 2)$ and $(3, 3)$ correspond, via set partitions,
to ovoids, perp sets and grids, respectively.

The Veldkamp sum operation on $V_{GQ(2,2)}$ described in the previous section
may now be defined purely in terms of sets: the Veldkamp sum of the two
set partitions $\{A|B\}$ and $\{C|D\}$ is given by $$
\{(A \cap C) \cup (B \cap D)|(A \cap D) \cup (B \cap C)\}
.$$

This identification extends to a set-theoretic description of the hyperplanes
of $L_3 \times$ GQ(2,\,2). The hyperplanes of this larger space may be put
into bijection with ordered quadruples of pairwise disjoint sets $(A, B, C, D)$
such that (a) no three of the sets are empty and (b) the union of the four
sets is $\{1, 2, 3, 4, 5, 6\}$. Such a quadruple corresponds to the hyperplane
given by the ordered triple of partitions $$
(\{A \cup B|C \cup D\}, \{A \cup C|B \cup D\}, \{A \cup D|B \cup C\})
.$$ Here, the leftmost component of the 
ordered triple describes the hyperplane of GQ(2,\,2) appearing in the uppermost
GQ(2,\,2)-quad of $L_3 \times$ GQ(2,\,2), and so on. For example, if
the sets $C$ and $D$ are empty, the top GQ(2,\,2)-quad will be deep
 and the other two will be identical to each other, being either singular, sub-quadrangular or ovoidal. 

The correspondence between the ordered quadruples and the hyperplanes is
four-to-one, because the quadruples $(A, B, C, D)$, $(B, A, D, C)$, 
$(C, D, A, B)$ and $(D, C, B, A)$ all index the same hyperplane. It follows
that acting by an element of the Klein four-group $V_4$ on an ordered quadruple
leaves the corresponding hyperplane invariant. The group $S_6 \times S_4$
acts on the quadruples, where $S_6$ acts diagonally on each of the set
partitions $A$, $B$, $C$ and $D$, and $S_4$ acts by place permutation. 
This induces an action of $S_6 \times S_4$ on the hyperplanes of 
$L_3 \times$ GQ(2,\,2), and since the action of $V_4 \leq S_4$ is trivial, 
this in turn induces an action of $S_6 \times (S_4/V_4) \cong S_6 \times S_3$ 
on the hyperplanes, thus recovering the full automorphism group of 
$L_3 \times$ GQ(2,\,2) in which $S_3$ acts by permuting the GQ(2,\,2)-quads.

This approach yields another way to deduce that the number of hyperplanes
of $L_3 \times$ GQ(2,\,2) is $2^{10} - 1$, as follows.  There are $4^6$ 
possible quadruples of pairwise disjoint sets $(A, B, C, D)$ whose union
is $\{1, 2, 3, 4, 5, 6\}$, and four of these quadruples have three empty
components. Since the correspondence between quadruples and hyperplanes is
four-to-one, the number of hyperplanes is $(4^6 - 4)/4$.

The correspondence described above induces a natural correspondence between
$S_6 \times S_4$-orbits (or $S_6 \times S_3$-orbits) of hyperplanes on the
one hand, and partitions of $6$ into two, three or four parts on the other.
There are eight such partitions; they are shown in Table 1, together
with their orbit sizes, stabilizers isomorphism types, stabilizer orders, 
and their name in the $H_1-H_8$ notation of \cite[Table 2]{sng}.

\begin{table}
\caption{A classification of geometric hyperplanes of $L_3 \times$ GQ(2,\,2).}
\vspace*{.3cm}
\centering
\begin{tabular}{|l|l|r|l|r|} \hline
Name & Partition & Orbit size & Stabilizer & Order \\ \hline
$H_1$ & $(3,3)$ &     30 &        $(S_3 \wr {\Bbb Z}_2) \times S_2$ & 144 \\
$H_2$ & $(4,2)$ &     45 &        $S_4 \times S_2 \times S_2$      & 96 \\
$H_3$ & $(5,1)$ &     18 &        $S_5 \times S_2$  & 240 \\
$H_4$ & $(2,2,1,1)$ & 270 &       $S_2 \times S_2 \times S_2 \times S_2$ & 16 \\
$H_5$ & $(2,2,2)$ &   90 &        $S_2 \times S_2 \times S_2 \times S_3$ & 48 \\
$H_6$ & $(3,1,1,1)$ & 120 &       $S_3 \times S_3$           & 36 \\
$H_7$ & $(3,2,1)$ &   360 &       $S_3 \times S_2$           & 12 \\
$H_8$ & $(4,1,1)$ &  90 &         $S_4 \times S_2$           & 48 \\
\hline
\end{tabular}
\end{table}

\section{Counting and classifying different types of Veldkamp lines}

The orbits of lines in the Veldkamp space $V$ may be enumerated using a
standard technique sometimes (inaccurately) known as Burnside's Lemma, which
proves the following.

Let $G$ be a finite group acting on a finite set $X$ with $t$ orbits, 
and for each $g \in G$, let $X^g$ denote the number of elements of $X$ fixed
by $g$.  Then we have $\displaystyle{
t = \frac{1}{|G|} \sum_{g \in G} |X^g|
.}$  Furthermore, if ${\mathcal C}$ is a set of
conjugacy class representatives of $G$, then we have $$
t = \frac{1}{|G|} \sum_{g \in {\mathcal C}} |{\mathcal C}| |X^g|
.$$

Using this technique, we can recover known results about orbits of lines 
under the action of the automorphism group $S_6$ of GQ(2,\,2):
there are $3$ orbits of hyperplanes (Veldkamp points)
and $5$ orbits of Veldkamp lines.
We can also recover the result the Veldkamp space $V$ has $8$ orbits of
hyperplanes under the automorphism group $S_6 \times S_3$.

The same idea can be adapted to count the orbits of Veldkamp lines of $V$.
The counting argument is more complicated than for the case of Veldkamp 
points, because it is possible for a line to be fixed by a group element 
$g$ without the three individual points being fixed.
There are three possibilities to consider, which we denote by (1), (2) and
(3) in Table 2.

\begin{itemize}
\item[(1)]{Every point of the Veldkamp line is fixed by $g$.  Such lines lie 
entirely within the fixed point space of $g$. Each number in the Fix(1)
column is the number of lines in a projective space PG$(d(g)-1, 2)$, for
a suitable integer $d(g)$ depending on the conjugacy class of $g$.}
\item[(2)]{One point of the Veldkamp line is fixed by $g$, and the other two
are exchanged.  To enumerate such lines, we take one point $x$ \emph{outside} 
the fixed point space of $g$.  The other two points are the point $g(x)$, 
and the point collinear with both of them (which is fixed by $g$).  We then 
divide by $2$ to correct for the overcount. 

Writing $d(g)$ as above, it follows
in each case that the entry in the Fix(2) column of $g$ is given by $$
{1 \over 2} \left( 2^{d(g^2)} - 2^{d(g)} \right)
.$$}
\item[(3)]{The element $g$ rotates the three points of the Veldkamp line in
a $3$-cycle. Each entry in the Fix(3) column is a number of the form
$(4^k - 1)/3$, and the enumeration of these cases is the most complicated. 
An ordered Veldkamp line can be thought of as a sequence of 30 binary digits. 
Typically, some even number, $2k$, of these bits can be 
chosen arbitrarily, provided that not all of them are zero, and then the rest
of the structure is forced. It is then necessary to divide by 3 to correct
for an overcount, by identifying an ordered Veldkamp line with each of
its cyclic shifts.}
\end{itemize}

\begin{table}[pth!]
\caption{Orbits of Veldkamp lines of $L_3 \times$ GQ(2,\,2).}
\vspace*{.3cm}
\centering
\begin{tabular}{|c|c|c|c|r|r|} \hline
Conjugacy class & Fix(1) & Fix(2) & Fix(3) & Size of class & Product \\ \hline 
id & 174251 & 0 & 0 & 1 & 174251 \\
$(1 2)$ & 10795 & 384 & 0 & 15 & 167685 \\
$(1 2)(3 4)$ & 651 & 480 & 0 & 45 & 50895 \\
$(1 2)(3 4)(5 6)$ & 651 & 480 & 0 & 15 & 16965 \\
$(1 2 3)$ & 651 & 0 & 5 & 40 & 26240 \\
$(1 2 3)(4 5 6)$ & 1 & 0 & 85 & 40 & 3440 \\
$(1 2 3 4)$ & 35 & 24 & 0 & 90 & 5310 \\
$(1 2 3 4)(5 6)$ & 35 & 24 & 0 & 90 & 5310 \\
$(1 2 3)(4 5)$ & 35 & 24 & 5 & 120 & 7680 \\
$(1 2 3 4 5)$ & 1 & 0 & 0 & 144 & 144 \\
$(1 2 3 4 5 6)$ & 1 & 0 & 5 & 120 & 720 \\
$(7 8)$ & 155 & 496 & 0 & 3 & 1953 \\
$(1 2)(7 8)$ & 155 & 496 & 0 & 45 & 29295 \\
$(1 2)(3 4)(7 8)$ & 155 & 496 & 0 & 135 & 87885 \\
$(1 2)(3 4)(5 6)(7 8)$ & 155 & 496 & 0 & 45 & 29295 \\
$(1 2 3)(7 8)$ & 7 & 28 & 1 & 120 & 4320 \\
$(1 2 3)(4 5 6)(7 8)$ & 0 & 1 & 5 & 120 & 720 \\
$(1 2 3 4)(7 8)$ & 7 & 28 & 0 & 270 & 9450 \\
$(1 2 3 4)(5 6)(7 8)$ & 7 & 28 & 0 & 270 & 9450 \\
$(1 2 3)(4 5)(7 8)$ & 7 & 28 & 1 & 360 & 12960 \\
$(1 2 3 4 5)(7 8)$ & 0 & 1 & 0 & 432 & 432 \\
$(1 2 3 4 5 6)(7 8)$ & 0 & 1 & 5 & 360 & 2160 \\
$(7 8 9)$ & 0 & 0 & 341 & 2 & 682 \\
$(1 2)(7 8 9)$ & 0 & 0 & 85 & 30 & 2550 \\
$(1 2)(3 4)(7 8 9)$ & 0 & 0 & 21 & 90 & 1890 \\
$(1 2)(3 4)(5 6)(7 8 9)$ & 0 & 0 & 21 & 30 & 630 \\
$(1 2 3)(7 8 9)$ & 1 & 0 & 85 & 80 & 6880 \\
$(1 2 3)(4 5 6)(7 8 9)$ & 35 & 0 & 21 & 80 & 4480 \\
$(1 2 3 4)(7 8 9)$ & 0 & 0 & 5 & 180 & 900 \\
$(1 2 3 4)(5 6)(7 8 9)$ & 0 & 0 & 5 & 180 & 900 \\
$(1 2 3)(4 5)(7 8 9)$ & 1 & 0 & 21 & 240 & 5280 \\
$(1 2 3 4 5)(7 8 9)$ & 0 & 0 & 1 & 288 & 288 \\
$(1 2 3 4 5 6)(7 8 9)$ & 1 & 6 & 5 & 240 & 2880 \\
\hline
& & & & & 673920 \\
\hline
\end{tabular}
\end{table}

We identify the group $S_6 \times S_3$ in the obvious way with the 
subgroup of $S_9$ fixing setwise each of the subsets 
$\{1, 2, 3, 4, 5, 6\}$ and $\{7, 8, 9\}$.  Since there are $11$ partitions of
$6$ and $3$ partitions of $3$, it follows that $S_6 \times S_3$ has $33$
conjugacy classes, and it is straightforward to find conjugacy class
representatives. Table 2 shows the calculation for the Veldkamp lines of 
$L_3 \times$ GQ(2,\,2).  The grand total of $$
673920 = |S_6 \times S_3| \times 156 = 720 \times 6 \times 156
$$ proves that there are 156 orbits of Veldkamp lines of the
near hexagon.

All 156 types are then listed in Table 3. Here, each type is characterized by its composition (columns 9 to 16) and the properties of the core $\mathcal{C}$ of the line, that is the set of points that are common to all the three geometric hyperplanes of a line of the given type. In particular, for each type (column 1) we list the number of points (column 2) and lines (column 3) of the core
as well as the distribution of the orders of its points. The last three columns show the intersection of $\mathcal{C}$ with each of the three GQ(2,\,2)-quads. Here, `g-perp' stands for a perp-set in a certain GQ(2,\,1) located in the particular GQ(2,\,2), and `unitr/tritr' abbreviates a unicentric/tricentric triad. If two or more types happen to possess the same string of parameters, the distinction between them is given by an explanatory remark/footnote.

\vspace*{-0.3cm}
\begin{table}[pth!]
\begin{center}
\caption{The types of Veldkamp lines of $L_3$ $\times$ GQ(2,\,2).} \vspace*{-0.1cm}
\resizebox{\columnwidth}{!}{%
{\begin{tabular}{||l|r|r|r|r|r|r|c|c|c|c|c|c|c|c|c|l|l|l||} \hline \hline
\multicolumn{1}{||c|}{} & \multicolumn{1}{|c|}{} & \multicolumn{1}{|c|}{}  &  \multicolumn{5}{|c|}{}
& \multicolumn{8}{|c|}{}                            & \multicolumn{1}{|c|}{} &\multicolumn{1}{|c|}{} &\multicolumn{1}{|c||}{}\\
\multicolumn{1}{||c|}{} & \multicolumn{1}{|c|}{} & \multicolumn{1}{|c|}{}  &  \multicolumn{5}{|c|}{$\#$ of Points of Order} & \multicolumn{8}{|c|}{Composition}
 & \multicolumn{1}{|c|}{} &\multicolumn{1}{|c|}{} &\multicolumn{1}{|c||}{}\\
 \cline{4-16}
Tp & Pt & Ln  &  0 &  1 &  2 &  3 &  4 & $H_1$ & $H_2$ & $H_3$ & $H_4$ & $H_5$ & $H_6$ & $H_7$ & $H_8$ & 1st & 2nd & 3rd  \\
\hline \hline
1  & 27  & 27   &  0 &  0 &  0 & 27 &  0 & 3     &  --   &  --   &  --   &  --   &  --   &  --   &  --   &  grid    &  grid   & grid \\
\hline
2 & 25  & 24   &  0 &  0 & 10 & 10 &  5 & 2     &  1    &  --   &  --   &  --   &  --   &  --   &  --   &  full    &  g-perp & g-perp\\
\hline
3  & 23  & 19   &  0 &  0 & 12 & 11 &  0 & 2     &  --   &  --   &  1    &  --   &  --   &  --   &  --   &  grid    &  g-perp & grid \\
\hline
4  & 21  & 20   &  0 &  0 & 6  & 12 &  3 & --    &  3    &  --   &  --   &  --   &  --   &  --   &  --   &  full    &  line   & line \\
5 & 21  & 18   &  0 &  6 &  0 & 12 &  3 & 1     &  1    &  1    &  --   &  --   &  --   &  --   &  --   &  full    &  unitr  & unitr\\
6 & 21  & 18   &  0 &  6 &  0 & 12 &  3 & --    &  3    &  --   &  --   &  --   &  --   &  --   &  --   &  full    &  tritr  & tritr\\
7 & 21  & 16   &  0 &  2 & 12 & 6  &  1 & 1     &  1    &  --   &  1    &  --   &  --   &  --   &  --   &  perp    &  grid   & g-perp\\
8  & 21  & 16   &  0 &  0 & 18 & 0  &  3 & --    &  3    &  --   &  --   &  --   &  --   &  --   &  --   &  perp    &  perp   & perp \\
\hline
9  & 19  & 15   &  0 &  0 & 12 & 7  &  0 & 1     &  --   &  --   &  2    &  --   &  --   &  --   &  --   &  grid    &  g-perp & g-perp \\
10  & 19  & 13   &  0 &  4 & 10 &  5 &  0 & 1    &  --   &  --   &  2    &  --   &  --   &  --   &  --   &  grid    &  g-perp & g-perp \\
11 & 19  & 12   &  0 &  6 &  9 & 4  &  0 & 1     &  1    &  --   &  --   &  --   &  --   &  1    &  --   &  perp    &  grid   & unitr\\
\hline
12 & 17  & 16   &  0 &  2 &  0 & 14 &  1 & --    &  1    &  2    &  --   &  --   &  --   &  --   &  --   &  full    &  point  & point\\
13 & 17  & 12   &  0 &  2 & 12 & 2  &  1 & --    &  1    &  --   &  2    &  --   &  --   &  --   &  --   &  perp    &  g-perp & g-perp\\
14  & 17  & 12   &  0 &  2 & 11 & 4  &  0 & --    &  1    &  --   &  2    &  --   &  --   &  --   &  --   &  grid    &  line   & g-perp \\
15  & 17  & 10   &  0 &  8 & 6  & 2  &  1 & 1     &  --   &  --   &  1    &  1    &  --   &  --   &  --   &  g-perp  &  g-perp & perp \\
16  & 17  & 10   &  1 &  4 & 10 & 2  &  0 & 1     &  --   &  --   &  1    &  --   &  --   &  1    &  --   &  grid    &  unitr  & g-perp \\
17 & 17  & 10   &  0 &  8 &  7 & 0  &  2 & --    &  2    &  --   &  --   &  1    &  --   &  --   &  --   &  perp    &  line   & perp\\
18  & 17  & 10   &  1 &  4 & 10 & 2  &  0 & --    &  1    &  --   &  2    &  --   &  --   &  --   &  --   &  grid    &  tritr  & g-perp \\
19  & 17  & 10   &  0 &  8 &  6 &  2 &  1 & --   &  1    &  --   &  2    &  --   &  --   &  --   &  --   &  perp    &  g-perp & g-perp\\
20 & 17  & 9    &  2 &  6 &  6 & 3  &  0 & 1     &  --   &  1    &  --   &  --   &  --   &  1    &  --   &  ovoid   &  unitr  & grid\\
21 & 17  & 9    &  0 &  8 &  8 & 1  &  0 & 1     &  --   &  --   &  1    &  --   &  1    &  --   &  --   &  perp    &  g-perp & g-perp\\
22 & 17  & 9    &  0 &  9 &  6 & 2  &  0 & --    &  2    &  --   &  --   &  --   &  1    &  --   &  --   &  perp    &  tritr   & perp\\
\hline
23 & 15  & 11   &  0 &  0 & 12 & 3  &  0 & --    &  --   &  --   &  3    &  --   &  --   &  --   &  --   &  g-perp$^{\hspace*{1.5mm}}$  &  g-perp & g-perp\\
24 & 15  & 9    &  0 &  6 & 6  & 3  &  0 & 1     &  --   &  --   &  --   &  --   &  --   &  2    &  --   &  unitr   &  grid   & unitr \\
25 & 15  & 9    &  0 &  6 & 6  & 3  &  0 & --    &  --   &  --   &  3    &  --   &  --   &  --   &  --   &  g-perp$^{1}$  &  g-perp & g-perp \\
26  & 15  &  9   &  0 &  6 &  6 &  3 &  0 & --   &  --   &  --   &  3    &  --   &  --   &  --   &  --   &  g-perp$^{1}$  &  g-perp & g-perp \\
27  & 15  & 8    &  2 &  4 & 7  & 2  &  0 & --    &  1    &  --   &  1    &  --   &  --   &  1    &  --   &  grid    &  tritr  & unitr \\
28  & 15  & 8    &  2 &  3 &  9 & 1  &  0 & --    &  1    &  --   &  1    &  --   &  --   &  1    &  --   &  line    &  grid   & unitr\\
29  & 15  & 8    &  2 &  4 & 7  & 2  &  0 & --    &  --   &  1    &  2    &  --   &  --   &  --   &  --   &  grid    &  unitr  & unitr \\
30  & 15  &  8   &  0 &  6 &  9 &  0 &  0 & --   &  --   &  --   &  3    &  --   &  --   &  --   &  --   &  g-perp  &  g-perp & g-perp \\
31  & 15  & 7    &  1 &  8 & 5  & 1  &  0 & 1     &  --   &  --   &  --   &  --   &  1    &  1    &  --   &  perp    &  g-perp & unitr \\
32  & 15  & 7    &  4 &  2 & 8  & 1  &  0 & 1     &  --   &  --   &  --   &  --   &  --   &  2    &  --   &  unitr   &  grid   & unitr \\
33  & 15  & 7    &  1 &  8 & 5  & 1  &  0 & --    &  1    &  --   &  1    &  --   &  --   &  1    &  --   &  perp    &  unitr  & g-perp \\
34  & 15  &  7   &  0 &  9 &  6 &  0 &  0 & --   &  --   &  --   &  3    &  --   &  --   &  --   &  --   &  g-perp  &  g-perp & g-perp \\
35  & 15  & 6    &  2 & 10 & 1  & 2  &  0 & 1     &  --   &  --   &  --   &  1    &  --   &  1    &  --   &  perp    &  unitr  & g-perp \\
36  & 15  & 6    &  3 &  6 & 6  & 0  &  0 & 1     &  --   &  --   &  --   &  --   &  --   &  2    &  --   &  ovoid   &  g-perp & g-perp \\
37  & 15  & 6    &  2 &  9 & 3  & 1  &  0 & --    &  1    &  1    &  --   &  --   &  --   &  1    &  --   &  ovoid   &  unitr  & perp \\
38  & 15  & 5    &  0 & 15 &  0 & 0  &  0 & --    &  --   &  3    &  --   &  --   &  --   &  --   &  --   &  ovoid   &  ovoid  & ovoid \\
\hline
39 & 13  & 8    &  0 &  4 &  8 & 0  &  1 & --    &  1    &  --   &  --   &  2    &  --   &  --   &  --   &  perp    &  line   & line\\
40  & 13  & 8    &  0 &  3 & 9  & 1  &  0 & --    &  1    &  --   &  --   &  --   &  --   &  2    &  --   &  line    &  grid   & point \\
41 & 13  & 8    &  0 &  4 & 7  & 2  &  0 & --    &  --   &  --   &  2    &  1    &  --   &  --   &  --   &  line    &  g-perp & g-perp \\
42  & 13  & 7    &  2 &  2 & 8  & 1  &  0 & --    &  --   &  1    &  1    &  --   &  --   &  1    &  --   &  grid    &  unitr  & point \\
43 & 13  & 6    &  0 &  9 &  3 & 1  &  0 & --    &  1    &  --   &  --   &  --   &  2    &  --   &  --   &  perp    &  tritr  & tritr\\
44  & 13  & 6    &  0 &  9 & 3  & 1  &  0 & --    &  1    &  --   &  --   &  --   &  2    &  --   &  --   &  perp    &  line   & line \\
45  & 13  & 6    &  4 &  0 & 9  & 0  &  0 & --    &  1    &  --   &  --   &  --   &  --   &  2    &  --   &  point   &  grid   & tritr \\
46  & 13  & 6    &  0 & 10 & 2  & 1  &  0 & --    &  1    &  --   &  --   &  --   &  --   &  2    &  --   &  perp    &  g-perp & point \\
47 & 13  & 6    &  0 &  9 &  3 & 1  &  0 & --    &  1    &  --   &  --   &  --   &  --   &  2    &  --   &  perp    &  unitr  & unitr\\
48  & 13  &  6   &  1 &  6 &  6 &  0 &  0 & --   &  --   &  --   &  2    &  --   &  1    &  --   &  --   &  tritr   &  g-perp & g-perp\\
49 & 13  & 6    &  0 &  8 & 5  & 0  &  0 & --    &  --   &  --   &  2    &  --   &  1    &  --   &  --   &  line    &  g-perp & g-perp \\
50  & 13  &  6   &  1 &  6 &  6 &  0 &  0 & --   &  --   &  --   &  2    &  --   &  --   &  1    &  --   &  g-perp  &  g-perp & unitr\\
\hline \hline
\end{tabular}%
} }
\end{center}
\end{table}
\addtocounter{table}{-1}
\vspace*{-0.3cm}
\begin{table}[pth!]
\begin{center}
\caption{(Continued.)} \vspace*{-0.1cm}
\resizebox{\columnwidth}{!}{%
{\begin{tabular}{||r|r|r|r|r|r|r|c|c|c|c|c|c|c|c|c|l|l|l||} \hline \hline
\multicolumn{1}{||c|}{} & \multicolumn{1}{|c|}{} & \multicolumn{1}{|c|}{}  &  \multicolumn{5}{|c|}{}
& \multicolumn{8}{|c|}{}                            & \multicolumn{1}{|c|}{} &\multicolumn{1}{|c|}{} &\multicolumn{1}{|c||}{}\\
\multicolumn{1}{||c|}{} & \multicolumn{1}{|c|}{} & \multicolumn{1}{|c|}{}  &  \multicolumn{5}{|c|}{$\#$ of Points of Order} & \multicolumn{8}{|c|}{Composition}
 & \multicolumn{1}{|c|}{} &\multicolumn{1}{|c|}{} &\multicolumn{1}{|c||}{}\\
 \cline{4-16}
Tp  & Pt  & Ln   &  0 &  1 &  2 &  3 &  4 & $H_1$ & $H_2$ & $H_3$ & $H_4$ & $H_5$ & $H_6$ & $H_7$ & $H_8$ & 1st      & 2nd     & 3rd  \\
\hline \hline
51  & 13  & 5    &  2 &  8 &  2 & 1  &  0 & --    &  1    &  --   &  --   &  1    &  1    &  --   &  --   &  perp    &  line   & tritr\\
52 & 13  & 5    &  2 &  8 &  2 & 1  &  0 & --    &  --   &  1    &  1    &  --   &  1    &  --   &  --   &  perp    &  unitr  & unitr\\
53 & 13  & 5    &  2 &  8 &  2 & 1  &  0 & --    &  --   &  --   &  2    &  1    &  --   &  --   &  --   &  tritr   &  g-perp & g-perp\\
54 & 13  & 5    &  0 & 11 & 2  & 0  &  0 & --    &  --   &  --   &  2    &  1    &  --   &  --   &  --   &  line    &  g-perp & g-perp \\
55  & 13  &  5   &  2 &  7 &  4 &  0 &  0 & --   &  --   &  --   &  2    &  --   &  1    &  --   &  --   &  tritr   &  g-perp & g-perp\\
56  & 13  &  5   &  2 &  8 &  2 &  1 &  0 & --   &  --   &  --   &  2    &  --   &  --   &  1    &  --   &  g-perp  &  g-perp & unitr\\
57 & 13  & 5    &  2 &  7 &  4 & 0  &  0 & --    &  --   &  --   &  2    &  --   &  --   &  1    &  --   &  unitr   &  g-perp & g-perp\\
58 & 13  & 4    &  4 &  8 &  0 & 0  &  1 & 1     &  --   &  --   &  --   &  1    &  --   &  --   &  1    &  perp    &  unitr  & unitr\\
59 & 13  & 4    &  4 &  8 &  0 & 0  &  1 & --    &  1    &  1    &  --   &  --   &  --   &  --   &  1    &  perp    &  ovoid  & point\\
60  & 13  & 4    &  4 &  8 &  0 & 0  &  1 & --    &  1    &  --   &  1    &  --   &  --   &  --   &  1    &  perp    &  unitr  & unitr\\
61  & 13  & 4    &  4 &  8 & 0  & 0  &  1 & --    &  1    &  --   &  --   &  2    &  --   &  --   &  --   &  perp    &  tritr  & tritr \\
62 & 13  & 4    &  4 &  7 & 1  & 1  &  0 & --    &  1    &  --   &  --   &  --   &  2    &  --   &  --   &  tritr   &  tritr  & perp\\
63  & 13  & 4    &  4 &  7 & 1  & 1  &  0 & --    &  1    &  --   &  --   &  --   &  --   &  2    &  --   &  line    &  g-perp & ovoid\\
64  & 13  &  4   &  4 &  7 &  1 &  1 &  0 & --   &  1    &  --   &  --   &  --   &  --   &  2    &  --   &  perp    &  unitr  & unitr\\
65 & 13  & 4    &  4 &  6 &  3 & 0  &  0 & --    &  1    &  --   &  --   &  --   &  --   &  2    &  --   &  tritr   &  g-perp & ovoid\\
66  & 13  & 4    &  4 &  8 &  0 & 0  &  1 & --    &  --   &  1    &  1    &  1    &  --   &  --   &  --   &  perp    &  unitr  & unitr \\
67 & 13  & 3    &  6 &  6 &  0 & 1  &  0 & 1     &  --   &  --   &  --   &  --   &  1    &  --   &  1    &  perp    &  unitr  & unitr\\
68  & 13  & 3    &  6 &  6 & 0  & 1  &  0 & 1     &  --   &  --   &  --   &  --   &  --   &  1    &  1    &  ovoid   &  g-perp & unitr \\
\hline
69  & 11  & 6    &  2 &  0 & 9  & 0  &  0 & --    &  --   &  1    &  --   &  --   &  --   &  2    &  --   &  grid    &  point  & point \\
70 & 11  & 5    &  0 &  7 & 4  & 0  &  0 & --    &  --   &  --   &  1    &  --   &  --   &  2    &  --   &  g-perp  &  g-perp & point \\
71 & 11  & 4    &  2 &  7 &  1 & 1  &  0 & --    &  --   &  1    &  --   &  1    &  --   &  1    &  --   &  perp    &  unitr  & point\\
72  & 11  & 4    &  2 &  7 & 1  & 1  &  0 & --    &  --   &  --   &  1    &  1    &  --   &  1    &  --   &  line    &  g-perp & unitr \\
73  & 11  &  4   &  2 &  6 &  3 &  0 &  0 & --   &  --   &  --   &  1    &  1    &  --   &  1    &  --   &  line    &  unitr  & g-perp\\
74 & 11  & 4    &  2 &  6 & 3  & 0  &  0 & --    &  --   &  --   &  1    &  --   &  1    &  1    &  --   &  unitr   &  tritr  & g-perp \\
75  & 11  &  4   &  2 &  6 &  3 &  0 &  0 & --   &  --   &  --   &  1    &  --   &  1    &  1    &  --   &  line    &  unitr  & g-perp\\
76 & 11  & 4    &  2 &  6 & 3  & 0  &  0 & --    &  --   &  --   &  1    &  --   &  --   &  2    &  --   &  g-perp$^{2}$  &  unitr  & unitr \\
77  & 11  &  4   &  2 &  6 &  3 &  0 &  0 & --   &  --   &  --   &  1    &  --   &  --   &  2    &  --   &  g-perp$^{2}$  &  unitr  & unitr \\
78  & 11  &  4   &  1 &  8 &  2 &  0 &  0 & --   &  --   &  --   &  1    &  --   &  --   &  2    &  --   &  point   &  g-perp & g-perp \\
79  & 11  & 3    &  4 &  6 & 0  & 1  &  0 & --    &  1    &  --   &  --   &  --   &  --   &  1    &  1    &  perp    &  point  & unitr \\
80 & 11  & 3    &  4 &  6 &  0 & 1  &  0 & --    &  --   &  1    &  --   &  --   &  1    &  1    &  --   &  perp    &  unitr  & point\\
81 & 11  & 3    &  2 &  9 & 0  & 0  &  0 & --    &  --   &  1    &  --   &  --   &  --   &  2    &  --   &  unitr   &  unitr  & ovoid\\
82  & 11  & 3    &  4 &  6 & 0  & 1  &  0 & --    &  --   &  --   &  2    &  --   &  --   &  --   &  1    &  unitr   &  g-perp & unitr \\
83  & 11  & 3    &  4 &  6 & 0  & 1  &  0 & --    &  --   &  --   &  1    &  1    &  --   &  1    &  --   &  tritr   &  unitr  & g-perp \\
84  & 11  &  3   &  4 &  5 &  2 &  0 &  0 & --   &  --   &  --   &  1    &  --   &  1    &  1    &  --   &  tritr   &  g-perp & unitr \\
85 & 11  & 3    &  3 &  7 &  1 & 0  &  0 & --    &  --   &  --   &  1    &  --   &  1    &  1    &  --   &  line    &  g-perp & unitr\\
86  & 11  &  3   &  4 &  6 &  0 &  1 &  0 & --   &  --   &  --   &  1    &  --   &  --   &  2    &  --   &  unitr$^{3}$   &  g-perp & unitr \\
87 & 11  & 3    &  4 &  6 & 0  & 1  &  0 & --    &  --   &  --   &  1    &  --   &  --   &  2    &  --   &  unitr$^{3}$   &  g-perp & unitr \\
88  & 11  &  3   &  4 &  5 &  2 &  0 &  0 & --   &  --   &  --   &  1    &  --   &  --   &  2    &  --   &  g-perp$^{4}$  &  unitr & unitr \\
89  & 11  &  3   &  4 &  5 &  2 &  0 &  0 & --   &  --   &  --   &  1    &  --   &  --   &  2    &  --   &  g-perp$^{4}$  &  unitr & unitr \\
90 & 11  & 2    &  6 &  4 &  1 & 0  &  0 & --    &  1    &  --   &  --   &  --   &  --   &  1    &  1    &  line    &  unitr  & ovoid\\
91  & 11  &  2   &  6 &  4 &  1 &  0 &  0 & --   &  --   &  --   &  2    &  --   &  --   &  --   &  1    &  unitr   &  g-perp & unitr \\
92  & 11  &  2   &  6 &  4 &  1 &  0 &  0 & --   &  --   &  --   &  1    &  1    &  --   &  1    &  --   &  tritr   &  unitr  & g-perp \\
93 & 11  & 2    &  6 &  4 & 1  & 0  &  0 & --    &  --   &  --   &  1    &  --   &  1    &  1    &  --   &  tritr   &  g-perp & unitr \\
94  & 11  &  2   &  6 &  4 &  1 &  0 &  0 & --   &  --   &  --   &  1    &  --   &  --   &  2    &  --   &  g-perp  &  unitr  & unitr \\
95 & 11  & 1    &  8 &  3 &  0 & 0  &  0 & --    &  --   &  2    &  --   &  --   &  --   &  --   &  1    &  ovoid   &  point  & ovoid\\
96  & 11  & 1    &  8 &  3 & 0  & 0  &  0 & --    &  --   &  1    &  --   &  --   &  --   &  2    &  --   &  unitr   &  unitr  & ovoid\\
97 & 11  & 0    & 11 &  0 &  0 & 0  &  0 & 1     &  --   &  --   &  --   &  --   &  --   &  --   &  2    &  unitr   &  unitr  & ovoid\\
98 & 11  & 0    & 11 &  0 & 0  & 0  &  0 & --    &  1    &  --   &  --   &  --   &  --   &  1    &  1    &  tritr   &  ovoid  & unitr\\\hline
99 & 9   & 6    &  0 &  0 &  9 & 0  &  0 & --    &  --   &  --   &  --   &  3    &  --   &  --   &  --   &  line    &  line   & line\\
100 & 9   & 4    &  0 &  8 &  0 & 0  &  1 & --    &  1    &  --   &  --   &  --   &  --   &  --   &  2    &  perp    &  point  & point\\
\hline \hline
\end{tabular}%
}}
\end{center}
\end{table}
\addtocounter{table}{-1}
\vspace*{-0.3cm}
\begin{table}[pth!]
\begin{center}
\caption{(Continued.)} \vspace*{-0.1cm}
\resizebox{\columnwidth}{!}{%
{\begin{tabular}{||l|r|r|r|r|r|r|c|c|c|c|c|c|c|c|c|l|l|l||} \hline \hline
\multicolumn{1}{||c|}{} & \multicolumn{1}{|c|}{} & \multicolumn{1}{|c|}{}  &  \multicolumn{5}{|c|}{}
& \multicolumn{8}{|c|}{}                            & \multicolumn{1}{|c|}{} &\multicolumn{1}{|c|}{} &\multicolumn{1}{|c||}{}\\
\multicolumn{1}{||c|}{} & \multicolumn{1}{|c|}{} & \multicolumn{1}{|c|}{}  &  \multicolumn{5}{|c|}{$\#$ of Points of Order} & \multicolumn{8}{|c|}{Composition}
 & \multicolumn{1}{|c|}{} &\multicolumn{1}{|c|}{} &\multicolumn{1}{|c||}{}\\
 \cline{4-16}
Tp  & Pt  & Ln   &  0 &  1 &  2 &  3 &  4 & $H_1$ & $H_2$ & $H_3$ & $H_4$ & $H_5$ & $H_6$ & $H_7$ & $H_8$ & 1st      & 2nd     & 3rd  \\
\hline \hline
101 & 9   & 3    &  2 &  6 &  0 & 1  &  0 & --    &  --   &  1    &  --   &  --   &  1    &  --   &  1    &  perp    &  point  & point\\
102 & 9   & 3    &  2 &  6 &  0 & 1  &  0 & --    &  --   &  --   &  1    &  --   &  --   &  1    &  1    &  point   &  g-perp & unitr\\
103 & 9   & 3    &  0 &  9 &  0 & 0  &  0 & --    &  --   &  --   &  --   &  3    &  --   &  --   &  --   &  line    &  line   & line\\
104 & 9   & 3    &  2 &  5 &  2 & 0  &  0 & --    &  --   &  --   &  --   &  2    &  1    &  --   &  --   &  line    &  tritr  & line\\
105 & 9   & 3    &  0 &  9 &  0 & 0  &  0 & --    &  --   &  --   &  --   &  1    &  2    &  --   &  --   &  line    &  line   & line\\
106 & 9   & 3    &  2 &  5 & 2  & 0  &  0 & --    &  --   &  --   &  --   &  1    &  --   &  2    &  --   &  tritr   &  g-perp & point \\
107  &  9  &  3   &  1 &  7 &  1 &  0 &  0 & --   &  --   &  --   &  --   &  1    &  --   &  2    &  --   &  point   &  g-perp & line \\
108 & 9   & 3    &  0 &  9 &  0 & 0  &  0 & --    &  --   &  --   &  --   &  --   &  3    &  --   &  --   &  tritr   &  tritr  & tritr\\
109  & 9   & 3    &  1 &  7 & 1  & 0  &  0 & --    &  --   &  --   &  --   &  --   &  1    &  2    &  --   &  point   &  g-perp & line \\
110 & 9   & 3    &  0 &  9 &  0 & 0  &  0 & --    &  --   &  --   &  --   &  --   &  --   &  3    &  --   &  unitr   &  unitr  & unitr\\
111  & 9   & 2    &  4 &  4 & 1  & 0  &  0 & --    &  --   &  --   &  1    &  --   &  1    &  --   &  1    &  line    &  unitr  & unitr \\
112  &  9  &  2   &  4 &  4 & 1 & 0 &  0 & --   &  --   &  --   &  1    &  --   &  --   &  1    &  1    &  g-perp   &  point  & unitr \\
113  &  9  &  2   &  4 &  4 &  1 &  0 &  0 & --   &  --   &  --   &  --   &  1    &  --   &  2    &  --   &  line    &  unitr$^{5}$  & unitr\\
114  & 9   & 2    &  4 &  4 &  1 & 0  &  0 & --    &  --   &  --   &  --   &  1    &  --   &  2    &  --   &  line    &  unitr$^{5}$  & unitr\\
115 & 9   & 2    &  4 &  4 &  1 & 0  &  0 & --    &  --   &  --   &  --   &  1    &  2    &  --   &  --   &  tritr   &  tritr  & line\\
116 & 9   & 2    &  3 &  6 &  0 & 0  &  0 & --    &  --   &  --   &  --   &  --   &  3    &  --   &  --   &  line    &  line   & tritr\\
117  & 9   & 2    &  4 &  4 & 1  & 0  &  0 & --    &  --   &  --   &  --   &  --   &  1    &  2    &  --   &  tritr   &  g-perp & point \\
118 & 9   & 2    &  3 &  6 &  0 & 0  &  0 & --    &  --   &  --   &  --   &  --   &  1    &  2    &  --   &  tritr   &  unitr  & unitr\\
119 & 9   & 2    &  4 &  4 & 1  & 0  &  0 & --    &  --   &  --   &  --   &  --   &  --   &  3    &  --   &  point$^{6}$   &  g-perp & unitr\\
120  & 9   &  2   &  4 &  4 &  1 &  0 &  0 & --   &  --   &  --   &  --   &  --   &  --   &  3    &  --   &  point$^{6}$   &  g-perp & unitr\\
121  & 9   & 1    &  6 &  3 & 0  & 0  &  0 & --    &  --   &  --   &  1    &  1    &  --   &  --   &  1    &  unitr   &  line   & unitr \\
122 & 9   & 1    &  6 &  3 &  0 & 0  &  0 & --    &  --   &  --   &  --   &  3    &  --   &  --   &  --   &  tritr   &  tritr  & line\\
123 & 9   & 1    &  6 &  3 &  0 & 0  &  0 & --    &  --   &  --   &  --   &  1    &  2    &  --   &  --   &  line    &  tritr  & tritr\\
124 & 9   & 1    &  6 &  3 & 0  & 0  &  0 & --    &  --   &  --   &  --   &  1    &  --   &  2    &  --   &  line    &  unitr  & unitr \\
125 & 9   & 1    &  6 &  3 &  0 & 0  &  0 & --    &  --   &  --   &  --   &  --   &  3    &  --   &  --   &  tritr   &  tritr  & tritr\\
126  & 9   & 1    &  6 &  3 & 0  & 0  &  0 & --    &  --   &  --   &  --   &  --   &  1    &  2    &  --   &  line    &  unitr  & unitr \\
127 & 9   & 1    &  6 &  3 & 0  & 0  &  0 & --    &  --   &  --   &  --   &  --   &  1    &  2    &  --   &  tritr   &  unitr  & unitr\\
128  &  9  &  1   &  6 &  3 &  0 &  0 &  0 & --   &  --   &  --   &  --   &  --   &  --   &  3    &  --   &  unitr   &  unitr  & unitr\\
129 & 9   & 0    &  9 &  0 &  0 & 0  &  0 & --    &  1    & --    &  --   &  --   &  --   &  --   &  2    &  tritr   &  point  & ovoid\\
130  & 9   & 0    &  9 &  0 & 0  & 0  &  0 & --    &  --   &  1    &  --   &  --   &  --   &  1    &  1    &  ovoid   &  unitr  & point \\
131 & 9   & 0    &  9 &  0 &  0 & 0  &  0 & --    &  --   &  --   &  1    &  1    &  --   &  --   &  1    &  tritr   &  unitr  & unitr\\
132 & 9   & 0    &  9 &  0 &  0 & 0  &  0 & --    &  --   &  --   &  1    &  --   &  1    &  --   &  1    &  tritr   &  unitr  & unitr\\
133  &  9  &  0   &  9 &  0 &  0 &  0 &  0 & --   &  --   &  --   &  1    &  --   &  --   &  1    &  1    &  unitr$^{7}$   &  unitr  &  unitr \\
134  &  9  &  0   &  9 &  0 &  0 &  0 &  0 & --   &  --   &  --   &  1    &  --   &  --   &  1    &  1    &  unitr$^{7}$   &  unitr  &  unitr \\
135 & 9   & 0    &  9 &  0 &  0 & 0  &  0 & --    &  --   &  --   &  --   &  2    &  1    &  --   &  --   &  tritr   &  tritr  & tritr\\
136 & 9   & 0    &  9 &  0 & 0  & 0  &  0 & --    &  --   &  --   &  --   &  1    &  --   &  2    &  --   &  tritr   &  unitr  & unitr \\
137  &  9  &  0   &  9 &  0 &  0 &  0 &  0 & --   &  --   &  --   &  --   &  --   &  1    &  2    &  --   &  tritr   &  unitr  & unitr \\
138  &  9  &  0   &  9 &  0 &  0 &  0 &  0 & --   &  --   &  --   &  --   &  --   &  --   &  3    &  --   &  unitr   &  unitr  & unitr\\
\hline
139 & 7   & 2    &  2 &  4 &  1 & 0  &  0 & --    &  --   &  --   &  --   &  1    &  --   &  1    &  1    &  point   &  unitr  & line\\
140  &  7  &  2   &  2 &  4 &  1 &  0 &  0 & --   &  --   &  --   &  --   &  --   &  --   &  2    &  1    &  point   &  g-perp & point \\
141  & 7   & 1    &  4 &  3 & 0  & 0  &  0 & --    &  --   &  1    &  --   &  --   &  --   &  --   &  2    &  ovoid   &  point  & point \\
142  & 7   & 1    &  4 &  3 & 0  & 0  &  0 & --    &  --   &  --   &  --   &  --   &  1    &  1    &  1    &  line    &  unitr  & point \\
143 & 7   & 1    &  4 &  3 & 0  & 0  &  0 & --    &  --   &  --   &  --   &  --   &  --   &  2    &  1    &  unitr$^{8}$   &  unitr  & point \\
144  &  7  &  1   &  4 &  3 &  0 &  0 &  0 & --   &  --   &  --   &  --   &  --   &  --   &  2    &  1    &  point   &  unitr$^{8}$  & unitr \\
145  & 7   & 0    &  7 &  0 & 0  & 0  &  0 & --    &  --   &  --   &  1    &  --   &  --   &  --   &  2    &  unitr   &  unitr  & point \\
146 & 7   & 0    &  7 &  0 &  0 & 0  &  0 & --    &  --   &  --   &  --   &  1    &  --   &  1    &  1    &  tritr   &  point  & unitr\\
147 & 7   & 0    &  7 &  0 &  0 & 0  &  0 & --    &  --   &  --   &  --   &  --   &  1    &  1    &  1    &  tritr   &  point$^{9}$  & unitr\\
148  &  7  &  0   &  7 &  0 &  0 &  0 &  0 & --   &  --   &  --   &  --   &  --   &  1    &  1    &  1    &  tritr   &  point$^{9}$  & unitr\\
149  &  7  &  0   &  7 &  0 &  0 &  0 &  0 & --   &  --   &  --   &  --   &  --   &  --   &  2    &  1    &  point$^{10}$   &  unitr  & unitr \\
150  &  7  &  0   &  7 &  0 &  0 &  0 &  0 & --   &  --   &  --   &  --   &  --   &  --   &  2    &  1    &  point$^{10}$   &  unitr$^{11}$  & unitr \\
\hline \hline
\end{tabular}%
}}
\end{center}
\end{table}
\addtocounter{table}{-1}
\begin{table}[pth!]
\begin{center}
\caption{(Continued.)} \vspace*{-0.1cm}
\resizebox{\columnwidth}{!}{%
{\begin{tabular}{||l|r|r|r|r|r|r|c|c|c|c|c|c|c|c|c|l|l|l||} \hline \hline
\multicolumn{1}{||c|}{} & \multicolumn{1}{|c|}{} & \multicolumn{1}{|c|}{}  &  \multicolumn{5}{|c|}{}
& \multicolumn{8}{|c|}{}                            & \multicolumn{1}{|c|}{} &\multicolumn{1}{|c|}{} &\multicolumn{1}{|c||}{}\\
\multicolumn{1}{||c|}{} & \multicolumn{1}{|c|}{} & \multicolumn{1}{|c|}{}  &  \multicolumn{5}{|c|}{$\#$ of Points of Order} & \multicolumn{8}{|c|}{Composition}
 & \multicolumn{1}{|c|}{} &\multicolumn{1}{|c|}{} &\multicolumn{1}{|c||}{}\\
 \cline{4-16}
Tp  & Pt  & Ln   &  0 &  1 &  2 &  3 &  4 & $H_1$ & $H_2$ & $H_3$ & $H_4$ & $H_5$ & $H_6$ & $H_7$ & $H_8$ & 1st      & 2nd     & 3rd  \\
\hline \hline
151  &  7  &  0   &  7 &  0 &  0 &  0 &  0 & --   &  --   &  --   &  --   &  --   &  --   &  2    &  1    &  point$^{10}$   &  unitr$^{11}$  & unitr \\
\hline
152  & 5   & 1    &  2 &  3 & 0  & 0  &  0 & --    &  --   &  --   &  --   &  1    &  --   &  --   &  2    &  line    &  point  & point \\
153 & 5   & 0    &  5 &  0 &  0 & 0  &  0 & --    &  --   &  --   &  --   &  --   &  1    &  --   &  2    &  tritr   &  point  & point\\
154 & 5   & 0    &  5 &  0 &  0 & 0  &  0 & --    &  --   &  --   &  --   &  --   &  --   &  1    &  2    &  unitr   &  point  & point\\
\hline
155 & 3   & 1    &  0 &  3 &  0 & 0  &  0 & --    &  --   &  --   &  --   &  --   &  --   &  --   &  3    &  point   &  point  & point\\
156 & 3   & 0    &  3 &  0 &  0 & 0  &  0 & --    &  --   &  --   &  --   &  --   &  --   &  --   &  3    &  point   &  point  & point\\
\hline \hline
\end{tabular}%
}}
\end{center}
\end{table}
\newpage

\noindent
Explanatory remarks:

$^{1}$Two (25) or no two (26) of the g-perps are such that their centers are joined by a type-one line.

$^{2}$The center of the g-perp does (77) or does not (76) lie on the type-one line passing through the center of
one of the two unicentric triads.

$^{3}$The centers of the two unicentric triads are (86) or are not (87) joined by a type-one line.

$^{4}$One line (88) or no line (89) of the g-perp is incident with the type-one line passing through the center of
one of the two unicentric triads.

$^{5}$The five type-one lines through the points of the two triads do (114) or do not (113) cut a doily-quad
in an ovoid.

$^{6}$One line (120) or no line (119) of type-two through the point is incident with the type-one line through the
center of the g-perp.

$^{7}$One (133) or none (134) of the unicentric triads is such that the type-one lines through two of its points
pass through the centers of the other two triads.

$^{8}$The centers of the two unicentric triads are (143) or are not (144) joined by a type-one line.

$^{9}$The point does (147) or does not (148) lie on the type-one line passing through a center of the tricentric triad.

$^{10}$The point does (149) or does not (150 and 151) lie on the type-one line passing through the center of one
of the two unicentric triads.

$^{11}$The centers of the two unicentric triads do (150) or do not (151) belong to the same grid-quad.

\section*{Acknowledgment}
This work already started in  2009, when the second author was a fellow of the Cooperation Group
``Finite Projective Ring Geometries: An Intriguing Emerging Link Between Quantum Information Theory, Black-Hole Physics and Chemistry of Coupling''  at the Center for Interdisciplinary Research (ZiF) of the University of Bielefeld, Bielefeld, Germany.
It was also supported in part by the VEGA Grant Agency, grant No. 2/0003/13.

\end{document}